\documentclass[10pt]{iopart}

\usepackage{iopams}  
\usepackage{graphics}
 
\usepackage{amssymb}
\usepackage{url}
\newtheorem{theorem}{Theorem}[section]

\newtheorem{lemma}[theorem]{Lemma}

\newtheorem{remark}[theorem]{Remark}
\newcommand{\R}{\mathbb{R}}
\renewcommand{\theequation}{\thesection.\arabic{equation}}
\providecommand{\abs}[1]{\left\vert#1\right\vert}
\providecommand{\norm}[1]{\left\Vert#1\right\Vert}
  
\renewcommand{\d}{\textrm{d}}
\renewcommand{\leq}{\leqslant}
\renewcommand{\geq}{\geqslant}
\def\epsilon{\varepsilon}
\def\phi {\varphi}

\begin{document}

\title[Stability of the determination of a coefficient for the wave equation]{Stability of the determination of a coefficient for the wave equation in an infinite wave guide}
\author{ Yavar Kian}
\address{ CPT, UMR CNRS 7332, Universit\'e d'Aix-Marseille, Universit\'e du Sud-Toulon-Var, CNRS-Luminy, 13288 Marseille, France.}
\ead{ yavar.kian@univ-amu.fr}
\begin{abstract}
We consider the stability in the inverse problem consisting in the determination of  an electric potential $q$, appearing in a Dirichlet initial-boundary value problem for  the wave equation $\partial_t^2u-\Delta u+q(x)u=0$ in an unbounded wave guide $\Omega=\omega\times\R$ with $\omega$ a bounded smooth domain of $\R^2$, from  boundary observations. The observation is given by the Dirichlet to Neumann map associated to a wave equation. We prove a H\"older stability estimate in the determination of $q$ from the  Dirichlet to Neumann map. Moreover, provided that the gap between two electric potentials rich its maximum in a fixed bounded subset of $\overline{\Omega}$, we extend this result  to the same inverse problem with  measurements on a bounded subset of the lateral boundary $(0,T)\times\partial\Omega$.
\end{abstract}

\maketitle
\section*{Introduction}

We consider the wave guide $\Omega=\omega\times\R$, where $\omega$ is a $\mathcal C^\infty$ bounded connected domain of  $\R^2$. We set  $\Sigma=(0,T)\times\partial\Omega$ and $Q=(0,T)\times\Omega$. Consider the  following initial-boundary
value problem (IBVP in short) for the wave equation
\begin{equation}\label{eq1}\left\{\begin{array}{ll}\partial_t^2u-\Delta u+q(x',x_3)u=0,\quad &t\in(0,T),\ x'\in\omega,\ x_3\in\R,\\  u(0,\cdot)=0,\quad \partial_tu(0,\cdot)=0,\quad &\textrm{in}\ \Omega,\\ u=f,\quad &\textrm{on}\ \Sigma.\end{array}\right.\end{equation}
Recall that $\partial\Omega=\partial\omega\times\R$. 
Since $\partial\Omega$ is not bounded, for all $s>0$ we give the following definition of the the space $H^s(\partial\Omega)$:
\[H^s(\partial\Omega)= H^s(\R_{{x}_3};L^2(\partial\omega))\cap L^2(\R_{{x}_3};H^s(\partial\omega)).\]
Then, we introduce the usual space
\[H^{r,s}((0,T)\times X)=H^r(0,T;L^2(X))\cap L^2(0,T;H^s(X))\]
where $X=\Omega$ or $X=\partial\Omega$.
Set the space
\[L=\left\{f\in H^{\frac{3}{2},\frac{3}{2}}(\Sigma): f_{\vert t=0}=0,\ \partial_tf,\partial_\tau f,\partial_{{x}_3}f\in L^2\left(\Sigma;\d\sigma(x)\frac{\d t}{t}\right)\right\}\]
with $\norm{}_L$ defined by
\[\norm{f}^2_L=\norm{f}^2_{H^{\frac{3}{2},\frac{3}{2}}(\Sigma)}+ \int_\Sigma\frac{\abs{\partial_tf}^2+\abs{\partial_\tau f}^2+\abs{\partial_{{x}_3}f}^2}{t}\d\sigma(x)\d t.\]
Here, we denote by $\partial_\tau$ a tangential derivative with respect to $\partial\omega$. We denote by $\nu$ the unit outward
normal vector to $\partial \Omega$. Notice that for $\nu_1$ the unit outward
normal vector to $\partial \omega$, we have $\nu=(\nu_1,0)$. We prove (see Theorem \ref{t2} in the appendix) that for $q\in L^\infty(\Omega)$ and $f\in L$ the IBVP (\ref{eq1}) has a unique solution
\[u_q\in \mathcal C([0,T];H^1(\Omega))\cap\mathcal C^1([0,T];L^2(\Omega))\]
such that $\partial_\nu u_q\in L^2(\Sigma)$. In addition, for any positive constant $M$, there exists a positive constant $C$ depending only of $\Omega$, $T$ and $M$, such that for all $q\in L^\infty(\Omega)$ with $\norm{q}_{L^\infty(\Omega)}\leq M$, the following estimate holds
\[\norm{u}_{\mathcal C([0,T];H^1(\Omega))}+\norm{u}_{\mathcal C^1([0,T];L^2(\Omega))}+\norm{\partial_\nu u_q}_{L^2(\Sigma)}\leq C\norm{f}_{L}.\]
In particular the following operator, usually called the Dirichlet to Neumann (DN map in short),
\[\begin{array}{rccl} \Lambda_q: & L& \to & L^2(\Sigma), \\
 \ \\ & f & \mapsto &\partial_\nu u_q \end{array}\]
is bounded. 

In the present paper, we consider the inverse problem which consists in determining the electric potential $q$ from the DN map $\Lambda_q$. We establish a stability estimate for this inverse problem. For $0<\alpha<1$ and $h\in\mathcal C(\overline{\Omega})$, we set
\[[h]_\alpha=\sup\left\{\frac{\abs{h(x)-h(y)}}{\abs{x-y}^\alpha}:\ x,y\in\overline{\Omega},\ x\neq y\right\}\]
and we consider the space 
\[\mathcal C_b^\alpha(\overline{\Omega})=\{h\in\mathcal C(\overline{\Omega})\cap L^\infty(\Omega):\ [h]_\alpha<\infty\}\]
with the norm
\[\norm{h}_{\mathcal C_b^\alpha(\overline{\Omega})}=\norm{h}_{L^\infty(\Omega)}+[h]_\alpha.\]
Our first main result can be stated as follows.

\begin{theorem} \label{t1} Let $M>0$, $0<\alpha<1$ and let $B_M$ be the ball centered at $0$ and of radius $M$ of $\mathcal C_b^\alpha(\overline{\Omega})$. Then, for $T>\textrm{Diam}(\omega)$ and  $q_1,q_2\in B_M$, 
we have
\begin{equation}\label{t1a}\norm{q_1-q_2}_{L^\infty(\Omega)}\leq C\norm{\Lambda_{q_1}-\Lambda_{q_2}}^{\frac{\min(2\alpha,1)\alpha}{3(2\alpha+2)\left(\min\left(4\alpha,2\right)+21\right)}}\end{equation}
with $C$ depending of $M$, $T$ and $\Omega$. Here $\norm{\Lambda_{q_1}-\Lambda_{q_2}}$ is the norm of $\Lambda_{q_1}-\Lambda_{q_2}$ with respect to $B\left(L,L^2(\Sigma)\right)$. \end{theorem}

Let us remark that in this result we consider the full DN map. This means that we determine the coefficient $q$ from measurements on the whole lateral boundary $\Sigma$ which is an unbounded set. This is due to the fact that we consider a large class of coefficients $q$ without  any restriction on their behavior outside a compact set (we only assume that the coefficients are uniformly bounded and H\"olderian). In order to extend this result to the determination of $q$ from measurements in a bounded subset of $\Sigma$, we need more informations about $q$. Namely, we need that the gap between  two coefficients $q_1,q_2$ reach its maximum in a fixed bounded subset of $\overline{\Omega}$. More precisely, let $R>0$ and consider the spaces $L_R$ which consists of functions $f\in L$ satisfying
\[f(t,x',x_3)=0,\quad t\in(0,T),\ x'\in\partial\omega,\ \abs{x_3}\geq R.\]
 Let us introduce the partial DN map defined by
 \[\begin{array}{rccl} \Lambda_q^R: & L_R& \to & L^2((0,T)\times\partial\omega\times (-R,R)), \\
 \ \\ & f & \mapsto &\partial_\nu {u_q}_{\vert (0,T)\times\partial\omega\times (-R,R)}. \end{array}\]
Our second result is the following.
 \begin{theorem} \label{tt1} Let $M>0$, $0<\alpha<1$ and let $B_M$ be the ball centered at $0$ and of radius $M$ of $\mathcal C_b^\alpha(\overline{\Omega})$. Let $T>\textrm{Diam}(\omega)$,  $q_1,q_2\in B_M$ and assume that there exists $r>0$ such that
\begin{equation}\label{tt1a}  \norm{q_1-q_2}_{L^\infty(\Omega)}=\norm{q_1-q_2}_{L^\infty(\omega\times (-r,r))}.\end{equation}
Then, for all $R>r$ we have
\begin{equation}\label{tt1b}\norm{q_1-q_2}_{L^\infty(\Omega)}\leq C\norm{\Lambda^R_{q_1}-\Lambda^R_{q_2}}^{\frac{\min(2\alpha,1)\alpha}{3(2\alpha+2)\left(\min\left(4\alpha,2\right)+21\right)}}\end{equation}
with $C$ depending of $M$, $T$, $\Omega$ and $R$. Here $\norm{\Lambda^R_{q_1}-\Lambda^R_{q_2}}$ is the norm of $\Lambda^R_{q_1}-\Lambda^R_{q_2}$ with respect to $B\left(L_R,L^2((0,T)\times\partial\omega\times (-R,R))\right)$. \end{theorem}

Clearly  condition (\ref{tt1a}) will be fulfilled if we assume that $q_1,q_2$ are compactly supported. Let us remark that this condition can also be fulfilled in more general cases. For example, consider the condition
\begin{equation}\label{per}v(x',x_3+2r)=v(x',x_3),\quad x'\in\overline{\omega},\ x_3\in\R.\end{equation}
Let $g:\R\to\R$ be a non negative continuous even function which is  decreasing in $(0,+\infty)$. Then, condition (\ref{tt1a}) will be fulfilled if we assume that $q_1,q_2$ are lying in the set
\[\hspace{-2cm}A_g=\{q:\ q(x',x_3)=g(x_3)v(x',x_3),\ v\in \mathcal C(\overline{\Omega})\cap L^\infty(\Omega),\ v\textrm{ satisfies (\ref{per})}\}.\]

In recent years the problem of recovering time-independent coefficients for hyperbolic equations in a bounded domain from boundary measurements  has attracted many attention.  In \cite{RS1}, the authors proved that the DN map determines uniquely the time-independent electric potential in a wave equation and \cite{RS2} has extended this result to the case of time-dependent potential. Isakov \cite{I} considered the determination of a coefficient of order zero and a damping coefficient. Note that all these results are concerned with measurements on  the whole boundary. The uniqueness by local DN map has been considered by \cite{E1} and \cite{E2}. The stability estimate in the case where the DN map is considered on the whole lateral boundary were treated by Stefanov and Uhlmann \cite{SU}. The uniqueness and H\"older stability estimate in a subdomain were established by Isakov and Sun \cite{IS} and, assuming that the coefficients are known in a neighborhood of the boundary, Bellassoued,  Choulli and  Yamamoto \cite{BCY} proved a log-type stability estimate in the case where the Neumann data are observed in an arbitrary subdomain of the boundary. In \cite{BJY1}, \cite{BJY2} and \cite{R} the authors established results with a finite number of data of DN map.

Let us also mention that the method using Carleman inequalities was first considered  by Bukhgeim and Klibanov \cite{BK}. For the application of Carleman estimate to the problem of recovering time-independent coefficients for hyperbolic equations we refer to \cite{B}, \cite{IY} and \cite{K}.

Let us observe that all these results are concerned with wave equations in a bounded domain. Several authors considered the problem of recovering time-independent coefficients in an unbounded domain from boundary measurements. Most of them considered the half space or the infinite slab. In \cite{R1}, Rakesh considered the problem of recovering the electric potential for the wave equation  in the half space from Neumann to Dirichlet  map. Applying   a   unique continuation result for the timelike Cauchy problem for the constant speed wave equation and the result of X-ray transform obtained by Hamaker, Smith, Solmon, Wagner in \cite{HSSW}, he proved a uniqueness result provided   that the electric potentials are  constant outside a compact set. In \cite{Na}, Nakamura extended this work to more general coefficients. In \cite{E3}, Eskin proved uniqueness modulo gauge invariance of magnetic and electric time-dependent potential with respect to the DN map  for the Schr\"odinger equation in a simply-connected bounded or unbounded domain. In \cite{Ik} and \cite{SW}, the authors considered the inverse problem of identifying an embedded object in an infinite slab.
In  \cite{LU}, the authors considered the problem of determining  coefficients for a stationary Schr\"odinger equation in an infinite slab. Assuming that the coefficients are compactly supported, they proved uniqueness  with respect to Dirichlet and Neumann data of the solution on parts of the boundary. This work was extended to the case of a magnetic stationary Schr\"odinger equation by \cite{KLU}.  In \cite{CS}, the authors considered the problem of determining the twisting for an elliptic equation in an infinite twisted wave guide. Assuming that  the first derivative of the twisting is sufficiently close to some \textit{a priori} fixed constant, they established a stability estimate of the twisting with respect to the DN map. To our best knowledge, with the one of \cite{CS}, this paper is the first where one establishes a stability estimate for the inverse problem of recovering a coefficient in an infinite domain with DN map without any assumption on the coefficient outside a compact set. 

The main ingredient in the proof of the stability estimates (\ref{t1a}) and (\ref{tt1b}) are  geometric optic solutions. The novelty in our approach comes from the fact that we take into account  the cylindrical form of the infinite wave guide and we use suitable  geometric optic solutions for this geometry.

This paper is organized as follows. In Section 1 we introduce the geometric optic solutions for our problem and, in a similar way to \cite{RS1} (see also Section 2.2.3 of \cite{Ch}), we prove existence of such solutions. Using these geometric optic solutions, in Section 2 we prove Theorem \ref{t1} and in Section 3 we prove Theorem \ref{tt1}. In the appendix, we treat the direct problem. We prove existence of solutions and we define the DN map. Let us remark that in the case of a bounded domain $\Omega$, applying some results of \cite{LLT}, \cite{BCY} have treated the direct problem. Since in this paper we consider an unbounded domain $\Omega$, it was necessary to treat this  problem.

 \section{Geometric optic solutions}
 
 The goal of this section is to construct geometric optic solutions  for the inverse problem. Let us recall that  every variable $x\in\Omega$ take the form $x=(x',x_3)$ with $x'\in\omega$ and $x_3\in\R$. Using this representation we can split the differential operator $\partial_t^2-\Delta$ defined on $Q$ into two differential operator $[\partial_t^2-\Delta_{x'}]+[-\partial_{{x}_3}^2]$ defined on $Q$. Keeping in mind this decomposition, we will construct geometric optic solutions $u\in H^2(Q)$ which are solutions in $Q$ of the equation
 $\partial_t^2u-\Delta u+qu=0$ and take the form 
  \[\hspace{-2,5cm}u(t,x',x_3)=\Phi(x'+t\theta)h(x_3)e^{\pm i\rho(x'\cdot\theta+t)}+\Psi^\pm(t,x',x_3;\rho),\  t\in(0,T),\ x'\in\omega,\ x_3\in \R\]
 with  $h\in\mathcal S(\R)$, $\Phi\in\mathcal C^\infty_0(\R^2)$, $\theta\in\mathbb S^1=\{y\in\R^2:\ \abs{y}=1\}$,  $\rho>0$ a parameter and $\Psi^\pm$ a remainder term that satisfies
 \[\norm{\Psi^\pm(.;\rho)}_{L^2(Q)}\leq\frac{C}{\rho}.\]
 Our result is the following.
 \begin{lemma}\label{l4} Let $q\in L^\infty(\Omega)$, $h\in\mathcal S(\R)$, $\Phi\in\mathcal C^\infty_0(\R^2)$, $\theta\in\mathbb S^1$, $\rho>0$ be arbitrary given.   Then the equation 
 \[\partial_t^2u-\Delta u+qu=0\]
 has solutions $u^\pm\in H^2(Q)$ of the form
 \[\hspace{-2,5cm}u^\pm(t,x',x_3)=\Phi(x'+t\theta)h(x_3)e^{\pm i\rho(x'\cdot\theta+t)}+\Psi^\pm(t,x',x_3;\rho),\ t\in(0,T),\ x'\in\omega,\ x_3\in \R.\]
 Here $\Psi^\pm$ satisfies
 \[\Psi^\pm(t,x;\rho)=0,\quad (t,x)\in\Sigma,\]
 \[\partial_t\Psi^+_{\vert t=0}=\Psi^+_{\vert t=0}=0,\]
  \[\partial_t\Psi^-_{\vert t=T}=\Psi^-_{\vert t=T}=0\]
and
\begin{equation}\label{l4a}\rho\norm{\Psi^\pm(.;\rho)}_{L^2(Q)}+\norm{\nabla_x\Psi^\pm(.;\rho)}_{L^2(Q)}\leq C\norm{h}_{H^2(\R)}\norm{\Phi}_{H^3(\R^2)},\end{equation}
where $C$ depends only on $T$, $\Omega$ and $M\geq \norm{q}_{L^\infty(\Omega)}$.
\end{lemma}
\textbf{Proof.}
We show existence of $u^+$. The existence of $u^-$ follows from similar arguments. Let $\Delta_{x'}$ be the Laplacian in $\omega$ and recall that \[[\Delta f](x',x_3)=[\Delta_{x'}f](x',x_3)+[\partial_{{x}_3}^2f](x',x_3),\quad f\in \mathcal C^2(\Omega),\ (x',x_3)\in\Omega.\]Notice that
\[\hspace{-2,5cm}(\partial_t^2-\Delta_{x'})\left[\Phi(x'+t\theta)h(x_3)e^{i\rho(x'\cdot\theta+t)}\right]=e^{i\rho(x'\cdot\theta+t)}\left[(\partial_t^2-\Delta_{x'})\Phi(x'+t\theta)h(x_3)\right]\]
and
\[\hspace{-2,5cm}-\partial_{{x}_3}^2\left[\Phi(x'+t\theta)h(x_3)e^{i\rho(x'\cdot\theta+t)}\right]=e^{i\rho(x'\cdot\theta+t)}\left[-\partial_{{x}_3}^2(\Phi(x'+t\theta)h(x_3))\right].\]
Therefore, we have
\[\hspace{-2,5cm}(\partial_t^2-\Delta+q)\left[\Phi(x'+t\theta)h(x_3)e^{i\rho(x'\cdot\theta+t)}\right]=e^{i\rho(x'\cdot\theta+t)}H(t,x',x_3)\]
with $H(t,x',x_3)=(\partial_t^2-\Delta+q)\Phi(x'+t\theta)h(x_3)$ and $\Psi^+$ must be the solution of
\begin{equation}\label{eq3}\hspace{-2cm}\left\{\begin{array}{ll}\partial_t^2\Psi-\Delta \Psi+q\Psi=e^{i\rho(x'\cdot\theta+t)}H(t,x',x_3),\quad &t\in(0,T),\ x'\in\omega,\ x_3\in\R,\\  \Psi(0,\cdot)=0,\quad \partial_t\Psi(0,\cdot)=0,\quad &\textrm{in}\ \Omega,\\ \Psi=0,\quad &\textrm{on}\ \Sigma.\end{array}\right.\end{equation}
Since $e^{i\rho(x'\cdot\theta+t)}H(t,x',x_3)\in L^2(Q)$, applying Theorem 8.1 in Chapter 3 of [LM1] (see also Remark \ref{r1}) we deduce the existence of $\Psi^+\in L^2(0,T;H^1_0(\Omega))\cap H^1(0,T;L^2(\Omega))$ solution of (\ref{eq3}). In addition, using the fact that $e^{i\rho(x'\cdot\theta+t)}H(t,x',x_3)\in H^1(0,T;L^2(\Omega))$, we can apply Theorem 2.1 in Chapter 5 of [LM2] (see Remark \ref{r2}) and prove that in fact  $\Psi^+\in H^2(Q)$.

\begin{remark}\label{r2} In order to apply Theorem 2.1 in Chapter 5 of [LM2] we combine the arguments introduced in Remark \ref{r1} with the fact that the operator $A=-\Delta+q$ with Dirichlet boundary condition is a selfadjoint operator  with domain $D(A)=H^2(\Omega)\cap H^1_0(\Omega)$. Then, we prove that Lemma 2.1 in Chapter 5 of [LM2] holds in our case and by the same way Theorem 2.1 in Chapter 5 of [LM2].\end{remark}
Now let 
\[W^+(t,x',x_3)=\int_0^t\Psi^+(s,x',x_3)\d s,\quad t\in(0,T),\ x'\in\omega,\ x_3\in \R.\]
Clearly $W^+$ is the solution of
\[\hspace{-2cm}\left\{\begin{array}{ll}\partial_t^2W-\Delta W+qW=\int_0^te^{i\rho(x'\cdot\theta+s)}H(s,x',x_3)\d s,\quad &t\in(0,T),\ x'\in\omega,\ x_3\in\R,\\  W(0,\cdot)=0,\quad \partial_tW(0,\cdot)=0,\quad &\textrm{in}\ \Omega,\\ W=0,\quad &\textrm{on}\ \Sigma.\end{array}\right.\]
From the energy estimate associated to the solution of this problem, we get
\[\norm{\Psi^+}_{L^2(Q)}=\norm{\partial_tW^+}_{L^2(Q)}\leq C\norm{\int_0^te^{i\rho(x'\cdot\theta+s)}H(s,x',x_3)\d s}_{L^2(Q)}.\]
Moreover, we have
\[\hspace{-2cm}\begin{array}{lll}\int_0^te^{i\rho(x'\cdot\theta+s)}H(s,x',x_3)\d s&=&\frac{1}{i\rho}\int_0^t\partial_se^{i\rho(x'\cdot\theta+s)}H(s,x',x_3)\d s\\
\ &=&-\frac{1}{i\rho}\int_0^te^{i\rho(x'\cdot\theta+s)}\partial_sH(s,x',x_3)\d s\\
\ &\ &+\frac{e^{i\rho(x'\cdot\theta+t)}H(t,x',x_3)-e^{i\rho(x'\cdot\theta)}H(0,x',x_3)}{i\rho}\end{array}\]
and it follows
\[\norm{\Psi^\pm(.;\rho)}_{L^2(Q)}\leq C\frac{\norm{h}_{H^2(\R)}\norm{\Phi}_{H^3(\R^2)}}{\rho}.\]
Combining this estimate with the energy estimate of (\ref{eq3}) we deduce (\ref{l4a}). \begin{flushright}
\rule{.05in}{.05in}
\end{flushright}

\section{Stability estimate}
The goal of this section is to prove Theorem \ref{t1}. 
Without lost of generality, we can assume that $0\in\omega$. From now on, we assume that $T>\textrm{Diam}(\omega)$, we fix $0< \epsilon<\min\left(1,\frac{T-\textrm{Diam}(\omega)}{3}\right)$ and  we set
\[\omega_\epsilon=\{ x\in\R^2\setminus\overline{\omega}:\ \textrm{dist}(x,\omega)<\epsilon\}.\]
We shall need a stability estimate for the problem of recovering a function from  X-ray transform.

\begin{lemma}\label{l5} Let $q_1,q_2\in L^\infty(\Omega)$, with $\norm{q_j}_{L^\infty(\Omega)}\leq M$, $j=1,2$, and let $q$ be equal to $q_1-q_2$ extended by $0$ outside of $\Omega$. Then, for all $\theta\in\mathbb S^1$ and $\Phi\in\mathcal C^\infty_0(\omega_\epsilon)$, $h\in\mathcal S(\R)$ we have
\begin{equation}\label{l5a}\begin{array}{l}\hspace{-2,5cm}\abs{\int_\R\int_{\R^2} \int_\R q(x',x_3)\Phi^2(x'+s\theta)h^2(x_3)\d s\d x'\d x_3}\\
\ \\
\hspace{-2,5cm}\leq C\left(\frac{\norm{h}^2_{H^2(\R)}\norm{\Phi}^2_{H^3(\R^2)}}{\rho}+\rho^2\norm{h}^2_{H^2(\R)}\norm{\Phi}^2_{H^3(\R^2)}\norm{\Lambda_{{q}_1}-\Lambda_{{q}_2}}\right),\quad \rho>1\end{array}\end{equation} with $C>0$ depending only of $\omega$, $M$ and $T$.\end{lemma}
\textbf{Proof.} In view of Lemma \ref{l4}, we can set 
 \[\hspace{-2cm}u_1(t,x',x_3;\theta,\rho)=\Phi(x'+t\theta)h(x_3)e^{i\rho(x'\cdot\theta+t)}+\Psi_1(t,x',x_3;\rho)\in H^2(Q),\]
  \[\hspace{-2cm}u_2(t,x',x_3;\theta,\rho)=\Phi(x'+t\theta)h(x_3)e^{-i\rho(x'\cdot\theta+t)}+\Psi_2(t,x',x_3;\rho)\in H^2(Q),\]
solutions of
 \[\partial_t^2u_1-\Delta u_1+q_1u_1=0,\quad \partial_t^2u_2-\Delta u_2+q_2u_2=0\]
with $\Psi_j\in H^2(Q)$, $j=1,2$, satisfying 
\begin{equation}\label{ll4a}\hspace{-2cm}\rho\norm{\Psi_j(.;\rho)}_{L^2(Q)}+\norm{\nabla_x\Psi_j(.;\rho)}_{L^2(Q)}\leq C\norm{h}_{H^2(\R)}\norm{\Phi}_{H^3(\R^2)},\ j=1,2,\end{equation}
 \[\Psi_j(.;\rho)_{\vert\Sigma}=0,
\quad j=1,2,\]
\[ \ \partial_t{\Psi_1}_{\vert t=0}={\Psi_1}_{\vert t=0}=0,
  \ \partial_t{\Psi_2}_{\vert t=T}={\Psi_2}_{\vert t=T}=0.\]
Since $T> \textrm{Diam}(\omega)+3\epsilon$,  we have
\[ \hspace{-2cm}\abs{t\theta +x'-y'}\geq t-\abs{x'-y'}\geq t-\textrm{Diam}(\omega)>2\epsilon,\quad x',y'\in\omega,\  t\geq T-\epsilon\]
and it follows
\[\hspace{-2cm}\{t\theta+x':\ \theta\in\mathbb S^1,\ x'\in\omega\}\subset \{y'\in\R^2:\ \textrm{dist}(y',\omega)>2\epsilon\}\subset \R^2\setminus\overline{\omega_\epsilon},\quad t\geq T-\epsilon.\]
Combining this with the  fact that $\textrm{supp}\Phi\subset \omega_\epsilon$,  we deduce
\begin{equation}\label{l5b}\partial_t^j\Phi(x'+t\theta)_{\vert t=s}=0,\quad x'\in\overline{\omega} ,\ s=T\ \textrm{or } s=0,\ j=0,1. \end{equation}
 Thus, we have 
 \begin{equation}\label{l5c}{u_1}_{\vert t=0}=\partial_t{u_1}_{\vert t=0}={u_2}_{\vert t=T}=\partial_t{u_2}_{\vert t=T}=0.\end{equation}
Let $f_\rho={u_1}_{\vert \Sigma}$   and notice that\[\hspace{-2cm}f_\rho(t,x',x_3)=\Phi(x'+t\theta)h(x_3)e^{i\rho(x'\cdot\theta+t)},\quad t\in(0,T),\  x'\in\partial\omega,\ x_3\in\R.\] In view of (\ref{l5b}) and Theorem 2.2 in Chapter 4 of [LM2], we have $f_\rho\in L$ and
\begin{equation}\label{l5d}\hspace{-2cm}\norm{f_\rho}_L\leq C\norm{\Phi(x'+t\theta)h(x_3)e^{i\rho(x'\cdot\theta+t)}}_{H^2(Q)}\leq C\rho^2\norm{h}_{H^2(\R)}\norm{\Phi}_{H^3(\R^2)}.\end{equation}
 Now let $v\in H^2(Q)$ be the solution of the IBVP
\[\left\{\begin{array}{ll}\partial_t^2v-\Delta v+q_2v=0,\quad &\textrm{in}\ Q,\\  v(0,\cdot)=0,\quad \partial_tv(0,\cdot)=0,\quad &\textrm{in}\ \Omega,\\ v=f_\rho,\quad &\textrm{on}\ \Sigma.\end{array}\right.\]
We set $u=v-u_1\in H^2(Q)$ and we have
\[\left\{\begin{array}{ll}\partial_t^2u-\Delta u+q_2u=qu_1,\quad &\textrm{in}\ Q,\\  u(0,\cdot)=0,\quad \partial_tu(0,\cdot)=0,\quad &\textrm{in}\ \Omega,\\ u=0,\quad &\textrm{on}\ \Sigma.\end{array}\right.\]
Applying (\ref{l5c}) and integrating by parts, we find
\[\int_Qqu_1u_2\d x\d t=-\int_\Sigma \partial_\nu uu_2\d\sigma( x)\d t.\]
Using the fact that $f_\rho\in L$, from this representation we get
\begin{equation}\label{l5e}\hspace{-2,5cm}\int_\Omega \int_0^T q(x',x_3)h^2(x_3)\Phi^2(x'+t\theta)\d t\d x'\d x_3=\int_QZ_\rho-\int_\Sigma(\Lambda_{{q}_2}f_\rho-\Lambda_{{q}_1}f_\rho)u_2\d\sigma(x)\d t\end{equation}
with \[\hspace{-2,5cm}\begin{array}{ll}Z_\rho(t,x',x_3)=-q(x',x_3)[&\Phi(x'+t\theta)h(x_3)\Psi_1(t,x',x_3)e^{-i\rho(x'\cdot\theta+t)}\\
\ &+\Phi(x'+t\theta)h(x_3)\Psi_2(t,x',x_3)e^{i\rho(x'\cdot\theta+t)}\\
\ &+\Psi_1(t,x',x_3)\Psi_2(t,x',x_3)].\end{array}\]
In view of (\ref{ll4a}), an application of the Cauchy-Schwarz inequality  yields
\[\int_Q\abs{Z_\rho}\d x\d t\leq \frac{2MC\norm{h}^2_{H^2(\R)}\norm{\Phi}^2_{H^3(\R^2)}}{\rho},\quad \rho>1\]
with $C$ depending of $\omega$, $T$.
From the fact that
\[\hspace{-2cm}u_2(t,x',x_3)=\Phi(x'+t\theta)h(x_3)e^{-i\rho(x'\cdot\theta+t)},\quad t\in(0,T),\  x'\in\partial\omega,\ x_3\in\R\]
and from (\ref{l5d}), we obtain
\[\hspace{-2cm}\begin{array}{ll}\abs{\int_\Sigma(\Lambda_{{q}_1}f_\rho-\Lambda_{{q}_2}f_\rho)u_2\d\sigma(x)\d t}&\leq \norm{\Lambda_{{q}_1}-\Lambda_{{q}_2}}\norm{f_\rho}_{L}\norm{\Phi(x'+t\theta)h(x_3)}_{L^2(\Sigma)}\\
\ &\leq C\rho^2\norm{\Lambda_{{q}_1}-\Lambda_{{q}_2}}\norm{h}^2_{H^2(\R)}\norm{\Phi}^2_{H^3(\R^2)}.\end{array}\]
Combining this estimates with (\ref{l5e}) and using the fact that $\textrm{supp}q\subset \Omega$ we get
\[\hspace{-2cm}\begin{array}{l}\abs{\int_\R\int_{\R^2} \int_0^T q(x',x_3)h^2(x_3)\Phi^2(x'+t\theta)\d t\d x'\d x_3}\\
\ \\
\leq C\left(\frac{\norm{h}^2_{H^2(\R)}\norm{\Phi}^2_{H^3(\R^2)}}{\rho}+\rho^2\norm{h}^2_{H^2(\R)}\norm{\Phi}^2_{H^3(\R^2)}\norm{\Lambda_{{q}_1}-\Lambda_{{q}_2}}\right),\quad \rho>1.\end{array}\]
Then,  using the fact that 
\[x'+t\theta\notin\overline{\omega_\epsilon},\quad x'\in\omega,\ t\geq T,\ \theta\in\mathbb S^1,\]
we find
\[\hspace{-2cm}\begin{array}{l}\abs{\int_\R\int_{\R^2} \int_0^{+\infty} q(x',x_3)h^2(x_3)\Phi^2(x'+t\theta)\d t\d x'\d x_3}\\
\ \\
\leq C\left(\frac{\norm{h}^2_{H^2(\R)}\norm{\Phi}^2_{H^3(\R^2)}}{\rho}+\rho^2\norm{h}^2_{H^2(\R)}\norm{\Phi}^2_{H^3(\R^2)}\norm{\Lambda_{{q}_1}-\Lambda_{{q}_2}}\right),\quad \rho>1.\end{array}\]
Repeating the above arguments with $v_j(t,x',x_3;\theta,\rho)=u_j(t,x',x_3; -\theta,\rho)$, $j=1,2$, we get
\[\hspace{-2cm}\begin{array}{l}\abs{\int_\R\int_{\R^2} \int_{-\infty}^{0} q(x',x_3)h^2(x_3)\Phi^2(x'+t\theta)\d t\d x'\d x_3}\\
\ \\
\leq C\left(\frac{\norm{h}^2_{H^2(\R)}\norm{\Phi}^2_{H^3(\R^2)}}{\rho}+\rho^2\norm{h}^2_{H^2(\R)}\norm{\Phi}^2_{H^3(\R^2)}\norm{\Lambda_{{q}_1}-\Lambda_{{q}_2}}\right),\quad \rho>1\end{array}\]
and we deduce (\ref{l5a}). \begin{flushright}
\rule{.05in}{.05in}
\end{flushright}

Let $\delta<\epsilon$. From now on, we will set the following. First, we fix $\phi\in\mathcal C^\infty_0(\R^2)$ real valued, supported on the  unit ball centered at $0$ and satisfying $\norm{\phi}_{L^2(\R^2)}=1$. We define
\[\Phi_\delta(x',y')=\delta^{-1}\phi\left(\frac{x'-y'}{\delta}\right),\quad x',y'\in\R^2.\]
We also define
\[h_\delta(x_3,y_3)=\delta^{-\frac{1}{2}}h\left(\frac{x_3-y_3}{\delta}\right),\quad x_3,y_3\in\R\]
with $h\in\mathcal C^\infty_0(\R)$ real valued,  supported on $[-1,1]$ and satisfying $\norm{h}_{L^2(\R)}=1$.
Finally, we set
\[\hspace{-2cm}R_\delta[q](y',y_3)=\int_\R\int_{\R^2} \Phi_\delta^2(x',y')h_\delta^2(x_3,y_3)q(x',x_3)\d x'\d x_3,\quad y'\in\R^2,\ y_3\in\R\]
and we introduce the  X-ray transform in $\R^2$ defined for all $f\in L^1(\R^2)$ by
\[[Xf](\theta,x')=\int_\R f(x'+t\theta)\d t,\quad x'\in\R^2,\ \theta\in\mathbb S^1.\]
We set also \[T\mathbb S^1=\{(x,\theta)\in\R^2\times\R^2:\ \theta\in\mathbb S^1,\ x\in\theta^{\bot}\}.\]
\begin{lemma}\label{l6} Let $M>0$, $0<\alpha<1$ and let $B_M$ be the ball centered at $0$ and of radius $M$ of $\mathcal C^\alpha_b(\overline{\Omega})$. Let $q_1,q_2\in B_M$ and let $q$ be equal to $q_1-q_2$ extended by $0$ outside of $\overline{\Omega}$. Then, for  $\delta^*=\frac{\epsilon}{4}$, we have
\begin{equation}\label{l6a}\hspace{-2,5cm}\norm{R_\delta[q]]}_{L^\infty\left(\R_{{y}_3}; L^2\left(\R^2_{y'}\right)\right)}\leq C\left(\frac{\delta^{-\frac{21}{2}}}{\rho}+\rho^2\delta^{-\frac{21}{2}}\norm{\Lambda_{{q}_1}-\Lambda_{{q}_2}}\right)^{\frac{1}{2}},\  0<\delta<\delta^*,\ \rho>1\end{equation}
with $C>0$ depending only of $\omega$, $M$ and $T$.\end{lemma}
\textbf{Proof.} Set 
\[\omega_1=\{x'\in\R^2:\ \frac{\epsilon}{4}<\textrm{dist}(x',\omega)<\frac{3\epsilon}{4}\},\]
\[\omega_2=\{x'\in\R^2:\ \textrm{dist}(x',\omega)<\frac{\epsilon}{4}\}\]
and let $\delta<\delta^*$. Then, for all $x'\in\R^2$ satisfying $\abs{x'-y'}\leq\delta$ with $y'\in \omega_1$ we have $x'\in\omega_\epsilon$.
Thus, $\Phi_\delta(.,y')\in\mathcal C^\infty_0(\omega_\epsilon)$ for all $y'\in\omega_1$.
In view of Lemma \ref{l5}, replacing $h$ and $\Phi$ by $h_\delta(.,y_3)$ and $\Phi_\delta(.,y')$ in (\ref{l5a}), for all $y'\in\omega_1$, $y_3\in\R$, $\rho>1$, we obtain
\[\hspace{-2,5cm}\begin{array}{l}\abs{\int_\R\int_{\R^2} \int_\R q(x',x_3)\Phi_\delta^2(x'+s\theta,y')h_\delta^2(x_3,y_3)\d s\d x'\d x_3}\\
\ \\
\leq C\left(\frac{\norm{h_\delta(.,y_3)}^2_{H^2(\R)}\norm{\Phi_\delta(.,y')}^2_{H^3(\R^2)}}{\rho}+\rho^2\norm{h_\delta(.,y_3)}^2_{H^2(\R)}\norm{\Phi_\delta(.,y')}^2_{H^3(\R^2)}\norm{\Lambda_{{q}_1}-\Lambda_{{q}_2}}\right).\end{array}\]
Applying the Fubini's theorem and the fact that $\Phi_\delta^2(x'+s\theta,y')=\Phi_\delta^2(x',y'-s\theta)$, we obtain
\[\hspace{-2,5cm}\begin{array}{l}\int_\R\int_{\R^2} \int_\R q(x',x_3)\Phi_\delta^2(x'+s\theta,y')h_\delta^2(x_3,y_3)\d s\d x'\d x_3\\
\ \\
= \int_\R\int_\R\int_{\R^2} q(x',x_3)\Phi_\delta^2(x',y'-s\theta)h_\delta^2(x_3,y_3)\d x'\d x_3\d s

=\int_\R R_\delta[q](y'-s\theta,y_3)\d s\\
\ \\

=X[R_\delta[q](.,y_3)](y',\theta)\end{array}\]
with $R_\delta[q](.,y_3)=y'\mapsto R_\delta[q](y',y_3)$.
Combining this representation with the previous estimate, for all $\rho>1$, we find
\[\hspace{-2,5cm}\begin{array}{l}\abs{X[R_\delta[q](.,y_3)](\theta,y')}\\
\ \\
\leq C\left(\frac{\norm{h_\delta(.,y_3)}^2_{H^2(\R)}\norm{\Phi_\delta(.,y')}^2_{H^3(\R^2)}}{\rho}+\rho^2\norm{h_\delta(.,y_3)}^2_{H^2(\R)}\norm{\Phi_\delta(.,y')}^2_{H^3(\R^2)}\norm{\Lambda_{{q}_1}-\Lambda_{{q}_2}}\right).\end{array}\]
One can easily check that
\[\hspace{-2cm}\norm{h_\delta(.,y_3)}^2_{H^2(\R)}\leq C\delta^{-4}, \ \norm{\Phi_\delta(.,y')}^2_{H^3(\R^3)}\leq C\delta^{-6},\quad y'\in\omega_1,\  y_3\in\R.\]
Thus,  we find
\begin{equation}\label{l6b}\hspace{-2,5cm}\abs{X[R_\delta[q](.,y_3)](\theta,y')}
\leq C\left(\frac{\delta^{-10}}{\rho}+\rho^2\delta^{-10}\norm{\Lambda_{{q}_1}-\Lambda_{{q}_2}}\right),\ \rho>1,\ y'\in\omega_1,\ y_3\in\R \end{equation}
with $C>0$ depending only of $\omega$, $M$ and $T$.
 From our choice for $\delta^*$ we get
\[\overline{\omega}+\{x'\in\R^2:\ \abs{x'}\leq\delta\}=\{x_1+x_2:\ x_1\in\overline{\omega},\ \abs{x_2}\leq\delta\}\subset \omega_2\]
and it follows
\begin{equation}\label{l6c}\textrm{supp}R_\delta[q](.,y_3)\subset\omega_2,\quad y_3\in\R.\end{equation}
Let $z'\in\omega_2$,  $\theta\in\mathbb S^1$ and consider the function $F$ defined on $\R^2$ by $F(x')=\textrm{dist}(x',\omega)$. One can easily check that
$F(\{z'+t\theta:\ t\in\R\})\supset[\frac{\epsilon}{4},+\infty)$ which implies that
$\{z'+t\theta:\ t\in\R\}\cap\omega_1\neq\varnothing$.
Therefore, for all $z'\in\omega_2$ and $\theta\in\mathbb S^1$ there exist $y'\in\omega_1$ and $t\in\R$ such that $z'=y'+t\theta$. Using the invariance property  of the X-ray transform we deduce that estimate (\ref{l6b}) holds for all $y'\in\omega_2$, $y_3\in\R$ and $\theta\in\mathbb S^1$. Now let $z'\in\R^2$ and $\theta\in\mathbb S^1$. If $\{z'+t\theta:\ t\in\R\}\cap\omega_2=\varnothing$, (\ref{l6c}) implies $X[R_\delta[q](.,,y_3)](z',\theta)=0$. On the other hand, if $\{z'+t\theta:\ t\in\R\}\cap\omega_2\neq\varnothing$ the invariance property of the X-ray transform implies that (\ref{l6b}) holds for $(z',\theta)$. Thus, (\ref{l6b}) holds for all $y'\in\R^2$, $\theta\in\mathbb S^1$ and $y_3\in\R$.
Let  $R>0$ be such that  $\overline{\omega_2}\subset B(0,R)=\{x'\in\R^2:\ \abs{x'}\leq R\}$.
For all $\theta\in\mathbb S^1$  we have
\[\abs{y'+t\theta}\geq\abs{y'}> R,\quad y'\in\theta^\bot\cap\R^2\setminus B(0,R),\ t\in\R\]
and (\ref{l6c}) implies
\[X[R_\delta[q](.,,y_3)](\theta,y')=0,\quad y'\in\theta^\bot\cap\R^2\setminus B(0,R).\]
It follows
\[\hspace{-2,5cm}\begin{array}{ll}\int_{T\mathbb S^1}\abs{X[R_\delta[q](.,y_3)]}^2\d y'\d\theta&=\int_{\mathbb S^1}\int_{\theta^\bot}\abs{X[R_\delta[q](.,y_3)](\theta,y')}^2\d y'\d \theta\\
\ \\
\ &=\int_{\mathbb S^1}\int_{\theta^\bot\cap B(0,R)}\abs{X[R_\delta[q](.,y_3)](\theta,y')}^2\d y'\d \theta,\quad y_3\in\R\end{array}\]
and (\ref{l6b}) implies
\[\int_{T\mathbb S^1}\abs{X[R_\delta[q](.,y_3)]}^2\d y'\d\theta\leq C\left(\frac{\delta^{-10}}{\rho}+\rho^2\delta^{-10}\norm{\Lambda_{{q}_1}-\Lambda_{{q}_2}}\right)^2.\]
In view of this estimate,  (\ref{l6c}) and the well known properties of the X-ray transform (see for example Theorem 5.1, p. 42 in \cite{N}), for all $y_3\in\R$, we obtain
\begin{equation}\label{l6i}\hspace{-2,5cm}\begin{array}{ll}\norm{R_\delta[q](.,y_3)}_{H^{-\frac{1}{2}}(\R^2)}&\leq C\norm{X[R_\delta[q](.,y_3)]}_{L^2(T\mathbb S^1)}\\
\ &\leq C\left(\frac{\delta^{-10}}{\rho}+\rho^2\delta^{-10}\norm{\Lambda_{{q}_1}-\Lambda_{{q}_2}}\right).\end{array}\end{equation}
Here $C$ depends only on $\omega$, $M$ and $T$. By interpolation we get,
\begin{equation}\label{l6d}\hspace{-2cm}\norm{R_\delta[q](.,y_3)}_{L^2(\R^2)}\leq C\norm{R_\delta[q](.,y_3)}_{H^{\frac{1}{2}}(\R^2)}^{\frac{1}{2}}\norm{R_\delta[q](.,y_3)}_{H^{-\frac{1}{2}}(\R^2)}^{\frac{1}{2}},\ y_3\in\R.\end{equation}
Consider  the following  estimate.
\begin{lemma} We have
\begin{equation}\label{l6e}\norm{R_\delta[q](.,y_3)}_{H^{\frac{1}{2}}(\R^2)}\leq C\delta^{-\frac{1}{2}},\ y_3\in\R\end{equation}
with $C$ depending only of $\omega$ and $M$.\end{lemma}
\textbf{Proof.} Note that
\[\int_{\R^2}\abs{R_\delta[q](y',y_3)}^2dy'=\int_{\R^2}\abs{\int_{\R^2}\Phi^2_\delta(x',y')G(x',y_3,\delta)dx'}^2dy'\]
with 
\[G(x',y_3,\delta)=\int_\R h^2_\delta(x_3,y_3)q(x',x_3)dx_3.\]
An application of the Cauchy-Schwarz inequality yields
\[\hspace{-2,5cm}\begin{array}{l}\int_{\R^2}\abs{R_\delta[q](y',y_3)}^2dy'\\
\ \\
\leq\int_{\R^2}\left[\left(\int_{\R^2}\Phi^2_\delta(x',y')dx'\right)\left(\int_{\R^2}\Phi^2_\delta(x',y')G^2(x',y_3,\delta)dx'\right)\right]dy'.\end{array}\]
Making the substitution $u= \frac{x'-y'}{\delta}$, we obtain
\[\hspace{-2,5cm}\int_{\R^2}\abs{R_\delta[q](y',y_3)}^2dy'\leq\int_{\R^2}\left[\left(\int_{\R^2}\phi^2(u)du\right)\left(\int_{\R^2}\phi^2(u)G^2(y'+\delta u,y_3,\delta)du\right)\right]dy'.\]
Then, applying  the Fubini's theorem and making the substitution $v=y'+\delta u$, we get
\begin{equation}\label{l6f}\hspace{-2,5cm}\begin{array}{ll}\int_{\R^2}\abs{R_\delta[q](y',y_3)}^2dy'&\leq\left(\int_{\R^2}\phi^2(u)du\right)^2\left(\int_{\R^2}G^2(v,y_3,\delta)dv\right)\\
\ \\
\ &\leq  \int_{\R^2}G^2(v,y_3,\delta)dv.\end{array}\end{equation}
Moreover, for all $v\in\R^2$, $y_3\in\R$, applying  the Cauchy-Schwarz inequality, we obtain
\[G^2(v,y_3,\delta)\leq \left(\int_\R h^2_\delta(x_3,y_3)dx_3\right)\left(\int_\R h^2_\delta(x_3,y_3)q^2(v,x_3)dx_3\right)\]
and, making the substitution $z_3=\frac{x_3-y_3}{\delta}$, we find
\[\hspace{-2,5cm}G^2(v,y_3,\delta)\leq \left(\int_\R h^2(z_3)dz_3\right)\left(\int_\R h^2(z_3)q^2(v,y_3+\delta z_3)dz_3\right)\leq \left(\int_\R h^2(z_3)dz_3\right)^24M^2.\]
Therefore, since $q=0$ outside of $\overline{\Omega}$, we have
\[\int_{\R^2}G^2(v,y_3,\delta)dv\leq \abs{\omega}\norm{G^2(.,y_3,\delta)}_{L^\infty(\R^2)}\leq 4\abs{\omega}M^2\]
and combining this estimate with (\ref{l6f}), we obtain 
\begin{equation}\label{l6g}\norm{R_\delta[q](.,y_3)}_{L^2(\R^2)}\leq C,\quad y_3\in\R\end{equation}
with $C$ depending of $\omega$ and $M$. Note that
\[\partial_{y'}R_\delta[q](y',y_3)=\int_{\R^2}\delta^{-3}G\left(\frac{x'-y'}{\delta}\right)G(x',y_3,\delta)dx'\]
with $G(u)=\partial_u\phi^2(u)$.
Therefore, repeating the above arguments, we get
\begin{equation}\label{l6h}\norm{R_\delta[q](.,y_3)}_{H^1(\R^2)}\leq C\delta^{-1},\quad y_3\in\R\end{equation}
with $C$ a constant depending of $\omega$ and $M$.
By interpolation, we obtain
\[\hspace{-1cm}\norm{R_\delta[q](.,y_3)}_{H^{\frac{1}{2}}(\R^2)}\leq C\norm{R_\delta[q](.,y_3)}^{\frac{1}{2}}_{L^2(\R^2)}\norm{R_\delta[q](.,y_3)}^{\frac{1}{2}}_{H^{1}(\R^2)},\quad y_3\in\R\]
and we deduce (\ref{l6e}) from (\ref{l6g}) and (\ref{l6h}).\begin{flushright}
\rule{.05in}{.05in}
\end{flushright}
In view of  estimates (\ref{l6i}), (\ref{l6d}) and (\ref{l6e}), we find
\[\norm{R_\delta[q](.,y_3)}_{L^2(\R^2)}\leq  C\left(\frac{\delta^{-\frac{21}{2}}}{\rho}+\rho^2\delta^{-\frac{21}{2}}\norm{\Lambda_{{q}_1}-\Lambda_{{q}_2}}\right)^{\frac{1}{2}},\ y_3\in\R\]
which implies (\ref{l6a}) since $C$ is independent of $y_3\in\R$. \begin{flushright}
\rule{.05in}{.05in}
\end{flushright}

\begin{lemma}\label{l7} Let $M>0$, $0<\alpha<1$ and let $B_M$ be the ball centered at $0$ and of radius $M$ of $\mathcal C^\alpha_b(\overline{\Omega})$. Let $q_1,q_2\in B_M$ and let $q$ be equal to $q_1-q_2$ extended by $0$ outside of $\overline{\Omega}$. Then, for  $\delta^*=\frac{\epsilon}{4}$ we have
\begin{equation}\label{l7a}\norm{R_\delta[q]-q}_{L^\infty\left(\R_{{y}_3}; L^2\left(\R^2_{y'}\right)\right)}\leq C\delta^{\tilde{\alpha}},\quad 0<\delta<\delta^*\end{equation}
with $C$ depending of $\omega$, $M$ and $\tilde{\alpha}=\min(\alpha,\frac{1}{2})$.
\end{lemma}
\textbf{Proof.}
Set 
\[S_\delta(y',y_3)=\int_{\R^2} \Phi_\delta^2(x',y')q(x',y_3)\d x',\quad y'\in\R^2,\ y_3\in\R.\]
Since $q\in \mathcal C^\alpha_b(\overline{\Omega}))$, one can easily check that for all $y_3\in\R$, $q(.,y_3)=y'\mapsto q(y',y_3)\in \mathcal C^\alpha(\overline{\omega})$. Thus, in view of Lemma 2.40 in [Ch], we have
\[\norm{S_\delta[q](.,y_3)-q(.,y_3)}_{ L^2(\R^2_{y'})}\leq C\delta^{\tilde{\alpha}}\norm{q(.,y_3)}_{ \mathcal C^\alpha\left(\overline{\omega}\right)}, y_3\in\R\]
with $C$ depending only of $\omega$, $M$ and $\alpha$.
It follows 
\begin{equation}\label{l7b}\hspace{-2,5cm}\norm{S_\delta[q]-q}_{L^\infty(\R_{{y}_3}; L^2(\R^2_{y'}))}\leq C\delta^{\tilde{\alpha}}\sup_{y_3\in\R}\norm{q(.,y_3)}_{ \mathcal C^\alpha\left(\overline{\omega}\right)}\leq C\delta^{\tilde{\alpha}}\norm{q}_{ \mathcal C_b^\alpha(\overline{\Omega})}\leq 2CM\delta^{\tilde{\alpha}}.\end{equation}
In view of this estimate, it only remains to prove
\begin{equation}\label{l7c}\norm{R_\delta[q]-S_\delta[q]}_{L^\infty(\R_{{y}_3}; L^2(\R^2_{y'}))}\leq C\delta^{\tilde{\alpha}}.\end{equation}
Notice that
\[\hspace{-2,5cm}R_\delta[q](y',y_3)-S_\delta[q](y',y_3)=\delta^{-1}\int_\R h^2\left(\frac{x_3-y_3}{\delta}\right) [S_\delta[q](y',x_3)-S_\delta[q](y',y_3)]\d x_3.\]
Making the substitution $u=\frac{x_3-y_3}{\delta}$ we find
\[\hspace{-2cm}R_\delta[q](y',y_3)-S_\delta[q](y',y_3)=\int_\R h^2(u)[S_\delta[q](y',y_3+\delta u)-S_\delta[q](y',y_3)]\d u.\]
In addition, for all $y'\in\R^2$, $y_3,u\in\R$, we have
\[\hspace{-2,5cm}\begin{array}{ll}\abs{S_\delta[q](y',y_3-\delta u)-S_\delta[q](y',y_3)}&\leq \int_{\R^2}\Phi_\delta^2(x',y')\abs{q(x',y_3+\delta u)-q(x',y_3)}dx'\\
\ \\
\ &\leq \int_{\omega}\Phi_\delta^2(x',y')\abs{q(x',y_3+\delta u)-q(x',y_3)}dx'\\
\ \\
\ &\leq \norm{q}_{\mathcal C^\alpha_b(\overline{\Omega})}(\delta \abs{u})^\alpha\int_{\omega}\Phi_\delta^2(x',y')dx'\\
\ \\
\ &\leq \norm{q}_{\mathcal C^\alpha_b(\overline{\Omega})}(\delta \abs{u})^\alpha\int_{\R^2}\Phi_\delta^2(x',y')dx'\\
\ \\
\ &\leq2M(\delta \abs{u})^\alpha.\end{array}\]
From this estimate we obtain
\[\norm{R_\delta[q]-S_\delta[q]}_{L^\infty(\R^3)}\leq 2M\delta^\alpha\int_\R h^2(u)\abs{u}^\alpha\d u=C\delta^\alpha\]
with $C$ depending only of $\omega$, $M$ and $\alpha$. Finally, since 
\[\hspace{-2cm}\textrm{supp}R_\delta[q](.,y_3)\cup \textrm{supp}S_\delta[q](.,y_3)\subset \overline{\omega}+\{x'\in\R^2:\ \abs{x'}\leq1\},\quad y_3\in\R,\]
for all $y_3\in\R$ we deduce
\[\hspace{-2cm}\begin{array}{ll}\norm{R_\delta[q](.,y_3)-S_\delta[q](.,y_3)}_{L^2(\R^2_{y'})}&\leq C\norm{R_\delta[q](.,y_3)-S_\delta[q](.,y_3)}_{L^\infty(\R^2_{y'})}\\
\ &\leq C\norm{R_\delta[q]-S_\delta[q]}_{L^\infty(\R^3)}\\
\ &\leq C\delta^\alpha.\end{array}\]
This last estimate implies (\ref{l7c}) and we deduce (\ref{l7a}).\begin{flushright}
\rule{.05in}{.05in}
\end{flushright}
\ \\
\textbf{Proof of Theorem \ref{t1} .} Let $q$ be equal to $q_1-q_2$ extended by $0$ outside of $\overline{\Omega}$. In view of (\ref{l6a}) and (\ref{l7a}), for $0<\delta<\delta^*$, we have
\begin{equation}\label{t1b}\hspace{-2,5 cm}\norm{q}_{L^\infty\left(\R_{{x}_3}; L^2\left(\R^2_{x'}\right)\right)}\leq C\left[\delta^{\tilde{\alpha}}+\left(\frac{\delta^{-\frac{21}{2}}}{\rho}+\rho^2\delta^{-\frac{21}{2}}\norm{\Lambda_{{q}_1}-\Lambda_{{q}_2}}\right)^{\frac{1}{2}}\right],\ \rho>1.\end{equation}
By interpolation, we obtain
\[\hspace{-2,5cm}\norm{q(.,y_3)}_{L^\infty(\omega)}\leq C\norm{q(.,y_3)}^{1-\mu}_{\mathcal C^\alpha(\overline{\omega})}\norm{q(.,y_3)}^\mu_{L^2(\omega)}\leq C(2M)^{1-\mu}\norm{q(.,y_3)}^\mu_{L^2(\omega)},\quad y_3\in\R\]
with $\mu=\frac{2\alpha}{2\alpha+2}$ and $C$ depending of $\omega$, $\alpha$. 
Combining this estimate with (\ref{t1b}), we obtain
\[\hspace{-2,5cm}\begin{array}{ll}\norm{q}_{L^\infty(\Omega)}&= \norm{q}_{L^\infty(\R_{{x}_3}; L^\infty(\omega))}\\
\ \\
\ &\leq C\left[\delta^{\tilde{\alpha}}+\left(\frac{\delta^{-\frac{21}{2}}}{\rho}+\rho^2\delta^{-\frac{21}{2}}\norm{\Lambda_{{q}_1}-\Lambda_{{q}_2}}\right)^{\frac{1}{2}}\right]^\mu,\ 0<\delta<\delta^*,\ \rho>1.\end{array}\]
This estimate can also be rewritten
\begin{equation}\label{t1c}\hspace{-2,5cm}\begin{array}{ll}\norm{q}_{L^\infty(\Omega)}\leq C2^{\frac{\mu}{2}}\left(\delta^{2\tilde{\alpha}}+\frac{\delta^{-\frac{21}{2}}}{\rho}+\rho^2\delta^{-\frac{21}{2}}\norm{\Lambda_{{q}_1}-\Lambda_{{q}_2}}\right)^{\frac{\mu}{2}},\ 0<\delta<\delta^*,\  \rho>1.\end{array}\end{equation}
Now let $\gamma=\norm{\Lambda_{{q}_1}-\Lambda_{{q}_2}}$ and set $\gamma^*={\delta^*}^{\frac{3(21+4\tilde{\alpha})}{2}}$. Then, for $0<\gamma<\gamma^*$ we can choose $\rho=\gamma^{-\frac{1}{3}}$, $\delta= \gamma^{\frac{2}{3(21+4\tilde{\alpha})}}$ and we obtain
\[\norm{q}_{L^\infty(\Omega)}\leq C\gamma^{\frac{2\tilde{\alpha}\alpha}{3(21+4\tilde{\alpha})(2\alpha+2)}},\quad \gamma<\gamma^*.\]
In addition, for $\gamma\geq\gamma^*$, we find
\[\norm{q}_{L^\infty(\Omega)}\leq\frac{2M}{{\gamma^*}^{\frac{2\tilde{\alpha}\alpha}{3(21+4\tilde{\alpha})(2\alpha+2)}}}\gamma^{\frac{2\tilde{\alpha}\alpha}{3(21+4\tilde{\alpha})(2\alpha+2)}}=C\gamma^{\frac{2\tilde{\alpha}\alpha}{3(21+4\tilde{\alpha})(2\alpha+2)}}.\]
Combining these two estimates we deduce (\ref{t1a}). \begin{flushright}
\rule{.05in}{.05in}
\end{flushright}

\section{Proof of Theorem \ref{tt1}}
 In this section we treat the special case introduced in Theorem \ref{tt1} where condition (\ref{tt1a}) is fulfilled. Like in the previous section, we assume that $0\in\omega$, $T>\textrm{Diam}(\omega)$ and we fix $0< \epsilon<\min\left(1,\frac{T-\textrm{Diam}(\omega)}{3}\right)$. We start with the following.
 
 \begin{lemma}\label{l8} Let $q_1,q_2\in L^\infty(\Omega)$, with $\norm{q_j}_{L^\infty(\Omega)}\leq M$, $j=1,2$, and let $q$ be equal to $q_1-q_2$ extended by $0$ outside of $\Omega$. Then, for all $\theta\in\mathbb S^1$ and $\Phi\in\mathcal C^\infty_0(\omega_\epsilon)$, $h\in\mathcal \mathcal C^\infty_0((-R,R))$ we have
\begin{equation}\label{l8a}\hspace{-2,5cm}\begin{array}{l}\abs{\int_\R\int_{\R^2} \int_\R q(x',x_3)\Phi^2(x'+s\theta)h^2(x_3)\d s\d x'\d x_3}\\
\ \\
\leq C\left(\frac{\norm{h}^2_{H^2(\R)}\norm{\Phi}^2_{H^3(\R^2)}}{\rho}+\rho^2\norm{h}^2_{H^2(\R)}\norm{\Phi}^2_{H^3(\R^2)}\norm{\Lambda^R_{{q}_1}-\Lambda^R_{{q}_2}}\right),\quad \rho>1\end{array}\end{equation}
with $C>0$ depending of $M$, $\omega$ and $T$.\end{lemma}
\textbf{Proof.} Repeating the arguments of Lemma \ref{l5}, we define
 \[\hspace{-2cm}u_1(t,x',x_3;\theta,\rho)=\Phi(x'+t\theta)h(x_3)e^{i\rho(x'\cdot\theta+t)}+\Psi_1(t,x',x_3;\rho)\in H^2(Q),\]
  \[\hspace{-2cm}u_2(t,x',x_3;\theta,\rho)=\Phi(x'+t\theta)h(x_3)e^{-i\rho(x'\cdot\theta+t)}+\Psi_2(t,x',x_3;\rho)\in H^2(Q),\]
solutions of
 \[\partial_t^2u_1-\Delta u_1+q_1u_1=0,\quad \partial_t^2u_2-\Delta u_2+q_2u_2=0\]
with $\Psi_j\in H^2(Q)$, $j=1,2$, satisfying 
\[\hspace{-2cm}\rho\norm{\Psi_j(.;\rho)}_{L^2(Q)}+\norm{\nabla_x\Psi_j(.;\rho)}_{L^2(Q)}\leq C\norm{h}_{H^2(\R)}\norm{\Phi}_{H^3(\R^2)},\ j=1,2,\]
 \[\Psi_j(.;\rho)_{\vert\Sigma}=0,
\quad j=1,2,\]
such that
 \[{u_1}_{\vert t=0}=\partial_t{u_1}_{\vert t=0}={u_2}_{\vert t=T}=\partial_t{u_2}_{\vert t=T}=0.\]
Let $f_j={u_j}_{\vert \Sigma}$, $j=1,2$,   and notice that\[\hspace{-2cm}f_j(t,x',x_3)=\Phi(x'+t\theta)h(x_3)e^{(-1)^ji\rho(x'\cdot\theta+t)},\quad t\in(0,T),\  x'\in\partial\omega,\ x_3\in\R,\] which implies that $f_j\in L_R $, $j=1,2$. 
 Now let $v\in H^2(Q)$ be the solution of the IBVP
\[\left\{\begin{array}{ll}\partial_t^2v-\Delta v+q_2v=0,\quad &t\in(0,T),\ x'\in\omega,\ x_3\in\R,\\  v(0,\cdot)=0,\quad \partial_tv(0,\cdot)=0,\quad &\textrm{in}\ \Omega,\\ v=f_1,\quad &\textrm{on}\ \Sigma.\end{array}\right.\]
Repeating the arguments of Lemma \ref{l5}, for $u=v-u_1$ we find \[\int_Qqu_1u_2\d x\d t=-\int_\Sigma \partial_\nu uu_2\d\sigma( x)\d t.\]
Since $f_2\in L_R$, we get
\[\int_\Sigma \partial_\nu uu_2\d\sigma( x)\d t=\int_0^T\int_{-R}^{R}\int_{\partial\omega}\partial_\nu uf_2\d\sigma(x')dx_3\d t\]
and since $\partial_\nu u=\Lambda_{{q}_2}f_1-\Lambda_{{q}_1}f_1$, with $f_1\in L_R $, we obtain
\[\hspace{-2,5cm}\int_Qqu_1u_2\d x\d t=-\int_\Sigma \partial_\nu uu_2\d\sigma( x)\d t=-\int_0^T\int_{-R}^{R}\int_{\partial\omega}[(\Lambda^R_{{q}_2}-\Lambda^R_{{q}_1})f_1]f_2\d\sigma(x')dx_3\d t.\]
Combining this representation with Lemma \ref{l5}, we prove easily estimate (\ref{l8a}). \begin{flushright}
\rule{.05in}{.05in}
\end{flushright}

\textbf{Proof of Theorem \ref{tt1}.} We consider the functions $\Phi_\delta(x',y')$, $h_\delta(x_3,y_3)$,
$R_\delta[q](y',y_3)$ introduced in the previous section and we extend  $q$  equal to $q_1-q_2$  by $0$ outside of $\overline{\Omega}$. Notice that for $\delta<\min\left(\frac{\epsilon}{4},R-r\right)$ and $y'\in\omega_1$ (with $\omega_1$ introduced in the proof of Lemma \ref{l6}), $y_3\in (-r,r)$ we have $\Phi_\delta(\cdot,y')\in\mathcal C^\infty_0(\omega_\epsilon)$ and $h_\delta(\cdot,y_3)\in\mathcal \mathcal C^\infty_0((-R,R))$. Therefore, combining Lemma \ref{l8} with Lemma \ref{l6}, we obtain the estimate
\begin{equation}\label{tt1c}\hspace{-2,5cm}\norm{R_\delta[q]}_{L^\infty(-r,r; L^2(\R^2_{y'}))}\leq C\left(\frac{\delta^{-\frac{21}{2}}}{\rho}+\rho^2\delta^{-\frac{21}{2}}\norm{\Lambda^R_{{q}_1}-\Lambda^R_{{q}_2}}\right)^{\frac{1}{2}},\  0<\delta<\delta^*,\ \rho>1\end{equation}
with $\delta^*=\min(\frac{\epsilon}{4},R-r)$. In addition, in view of Lemma \ref{l7}, we find
\begin{equation}\label{tt1d}\norm{R_\delta[q]-q}_{L^\infty(-r,r; L^2(\R^2_{y'}))}\leq C\delta^{\tilde{\alpha}},\quad 0<\delta<\delta^*.\end{equation}
Combining (\ref{tt1c}) and (\ref{tt1d}) with the arguments of Theorem \ref{t1}, we obtain easily
\begin{equation}\label{tt1e}\norm{q_1-q_2}_{L^\infty((-r,r)\times\omega)}\leq C\norm{\Lambda^R_{q_1}-\Lambda^R_{q_2}}^{\frac{\min(2\alpha,1)\alpha}{3(2\alpha+2)\left(\min\left(4\alpha,2\right)+21\right)}}\end{equation}
with $C$ depending of $R$, $M$, $\Omega$ and $T$. Then, condition (\ref{tt1a}) implies (\ref{tt1b}). \begin{flushright}
\rule{.05in}{.05in}
\end{flushright}

\section{Appendix}

In this appendix  we treat the direct problem. Our goal is to prove the following.

\begin{theorem} \label{t2} Let $q\in L^\infty(\Omega)$ and $f\in L$. Then problem (\ref{eq1}) admits a unique solution $u\in \mathcal C([0,T];H^1(\Omega))\cap\mathcal C^1([0,T];L^2(\Omega))$ such as $\partial_\nu u\in L^2(\Sigma)$. Moreover, this solution $u$ satisfies
\begin{equation}\label{t2a}\norm{u}_{\mathcal C([0,T];H^1(\Omega))}+\norm{u}_{\mathcal C^1([0,T];L^2(\Omega))}+\norm{\partial_\nu u}_{L^2(\Sigma)}\leq C\norm{f}_L.\end{equation}\end{theorem}

In the case of a bounded domain $\Omega$, applying Theorem 2.1 of \cite{LLT}, \cite{BCY} proved this result for $f\in H^1(\Sigma)$. Since $\Omega$ is an unbounded domain, we can not apply the analysis of \cite{LLT}. Nevertheless, we can  solve problem (\ref{eq1}) by the classical argument which comprises in transforming this problem into a problem with an inhomogeneous equation and homogeneous boundary conditions. For this propose, we first need to establish a result of lifting for  Sobolev spaces in a wave guide.

\subsection{Result of lifting for Sobolev spaces}

In this subsection we will show the following.
\begin{theorem} \label{t3} For all $f\in L$, there exists $w[f]\in H^{2,2}(Q)$ satisfying
\begin{equation}\label{t3a}\left\{\begin{array}{l} w[f](0,x)=\partial_tw[f](0,x)=0,\quad x\in\Omega,\\ \partial_\nu w[f]=0,\ w[f]=f,\ \quad\textrm{on }\Sigma\end{array}\right.\end{equation}
and
\begin{equation}\label{t3b}\norm{w[f]}_{H^{2,2}(Q)}\leq C\norm{f}_L.\end{equation}\end{theorem}

For this purpose, we will establish more general  result of lifting for Sobolev spaces. Repeating the arguments of pages 38-40 in Chapter 1 of \cite{LM1}, by the mean of local coordinates we can replace $\omega$ by $\R_+^2$, with $\R_+^2=\{(x_1,x_2)\in\R^3:\ x_1>0\}$, and $\partial\omega$ by $\R$. Using the fact that $\Omega=\omega\times\R$ and $\partial\Omega=\partial\omega\times\R$, with the same changes applied only to  $x'\in\overline{\omega}$ for any variable $x=(x',x_3)\in\overline{\Omega}$, we can replace $\Omega$ by $\R_+^3$, with $\R_+^3=\{(x_1,x_2,x_3)\in\R^3:\ x_1>0\}$, $\partial\Omega$ by $\R^2$, $\partial_\nu$ by $-\partial_{x_1}$ and, without lost of generality, we can assume $T=\infty$. Then, in our result we can replace  $H^{r,r}(Q)$ (respectively $H^{r}(\Omega)$ and $H^{r,r}(\Sigma)$) by $H^{r,r}((0,+\infty)\times\R_+^3)$ (respectively $H^{r}(\R_+^3)$ and $H^{r,r}((0,+\infty)\times\R^2)$). Let $K_1$ be the space of $(g_0,g_1,u_0,u_1)$ satisfying
\[g_j\in H^{r_j,r_j}((0,+\infty)\times\R^2),\ r_j=2-j-\frac{1}{2},\ j=0,1,\]
\[u_k\in H^{s_k}(\R_+^3),\ s_k=2-k-\frac{1}{2},\ k=0,1\]
and the compatibility conditions
\[ g_0(0,x_2,x_3)=u_0(0,x_2,x_3),\quad (x_2,x_3)\in \R^2,\]
\begin{equation}\label{com1}\hspace{-2,5cm}\begin{array}{l}\int_\R\int_\R\int_0^{+\infty}\abs{\partial_t^kg_j(.,x_2,x_3)_{\vert t=r}-\partial_{{x}_1}^ju_k(.,x_2,x_3)_{\vert x_1=r}}^2\frac{\d r}{r}\d x_2\d x_3,\\ j,k\in\mathbb N,\ j+k=1,\end{array}\end{equation}
\begin{equation}\label{com2}\hspace{-2,5cm}\begin{array}{l}\int_\R\int_\R\int_0^{+\infty}\abs{\partial_{({x}_2,{x}_3)}^\beta g_0(.,x_2,x_3)_{\vert t=r}-\partial_{({x}_2,{x}_3)}^{\beta}u_0(.,x_2,x_3)_{\vert {x}_1=r}}^2\frac{\d r}{r}\d x_2\d x_3,\\ \beta\in\mathbb N^2,\ \abs{\beta}=1.\end{array}\end{equation}
Conditions (\ref{com1}) and (\ref{com2}) are  global compatibility conditions (see subsection 2.4 in Chapter 4 of [LM2]). 
Let us also introduce the Hilbert space
\[ \hspace{-2cm}K_2=\{f\in H^{\frac{3}{2},\frac{3}{2}}((0,+\infty)\times\R^2):\ f_{\vert t=0}=0,\  t^{-\frac{1}{2}}\nabla_{(t,{x}_2,{x}_3)}f\in L^2((0,+\infty)\times \R^2)\}\]
with the norm
\[\hspace{-2cm}\norm{f}_{{K}_2}^2=\norm{f}_{H^{\frac{3}{2},\frac{3}{2}}((0,+\infty)\times\R\times \R)}^2+\norm{ t^{-\frac{1}{2}}\nabla_{(t,{x}_2,{x}_3)}f}_{L^2((0,+\infty)\times \R^2)}^2.\]
Using the above changes, we will deduce Theorem \ref{t3} from

\begin{lemma} \label{l1} The operator
\[ U: w\longmapsto(w_{\vert{x}_1=0}, \partial_{{x}_1}w_{\vert{x}_1=0}, w_{\vert t=0},\partial_tw_{\vert t=0})\]
is continuous and subjective from $H^{2,2}((0,+\infty)\times\R_+^3)$ to $K_1$.
In addition, for $f\in K_2$ there exists $w\in H^{2,2}((0,+\infty)\times\R_+^2\times \R)$ satisfying
\begin{equation}\label{l1f}(w_{\vert{x}_1=0}, \partial_{{x}_1}w_{\vert{x}_1=0}, w_{\vert t=0},\partial_tw_{\vert t=0})=(f,0,0,0)\end{equation}
with 
\begin{equation}\label{l1g}\norm{w}_{H^{2,2}((0,+\infty)\times\R_+^3)}\leq C\norm{f}_{K_2}.\end{equation}

\end{lemma}

\textbf{Proof.} According to Theorem 2.3 in Chapter 4 of \cite{LM2}, the operator $U$ is continuous and subjective. Therefore, it only remains to prove the last part of the lemma.   For this purpose, let $(g_0,g_1,u_0,u_1)=(f,0,0,0)$ and remark that $(g_0,g_1,u_0,u_1)\in K_1$.
In view of Theorem 4.2 in chapter 1 of \cite{LM1} and the first part of the proof of Theorem 2.1 in Chapter 4 of \cite{LM2},
there exists $v[f]\in H^{2,2}((0,+\infty)\times\R_+^3)$ such that 
\[(v[f]_{\vert {x}_1=0}, \partial_{{x}_1}v[f]_{\vert {x}_1=0})=(f,0)\]
and 
\begin{equation}\label{l1a}\norm{v[f]}_{H^{2,2}((0,+\infty)\times\R_+^3)}\leq C\norm{f}_{H^{\frac{3}{2},\frac{3}{2}}((0,+\infty)\times\R^2)}.\end{equation}
Consider $H^s_0(\R_+^3)$, $s>0$,  the closure of $\mathcal C^\infty_0(\R_+^3)$ in $H^s(\R_+^3)$ and let $H^{k+\frac{1}{2}}_{0,0}(\R_+^3)$, $k=0,1$, be the spaces
\[H^{1+\frac{1}{2}}_{0,0}(\R_+^3)=\{v\in H^{1+\frac{1}{2}}_{0}(\R_+^3):\ x_1^{-\frac{1}{2}}\partial_{{x}_j}v\in L^2(\R_+^3),\ j=1,2,3\},\]
\[H^{\frac{1}{2}}_{0,0}(\R_+^3)=\{v\in H^{\frac{1}{2}}_{0}(\R_+^3):\ x_1^{-\frac{1}{2}}v\in L^2(\R_+^3)\}.\]
According to Theorem 11.7 in Chapter 1 of \cite{LM1}, we have
\begin{equation}\label{l1b}[H^2_0(\R_+^3),L^2(\R_+^3)]_{\frac{j+\frac{1}{2}}{2}}=H^{2-j-\frac{1}{2}}_{0,0}(\R_+^3),\quad j=0,1,\end{equation}
where $[H^2_0(\R_+^3),L^2(\R_+^3)]_{\frac{j+\frac{1}{2}}{2}}$ is the interpolation space of order $\frac{j+\frac{1}{2}}{2}$ between $H^2_0(\R_+^3)$ and $L^2(\R_+^3)$.
In view of (\ref{l1b}) and   Theorem 3.2 in Chapter 1 of \cite{LM1}, the operator
\[u\longmapsto (u_{\vert t=0},\partial_tu_{\vert t=0})\]
is continuous and subjective from the Hilbert space\[L^2(\R^+_t;H^{2}_0(\R_+^3))\cap H^2(\R^+_t;L^2(\R_+^3)))\] to the Hilbert space
\[H^{\frac{3}{2}}_{0,0}(\R^3_+)\times H^{\frac{1}{2}}_{0,0}(\R^3_+).\]
Thus, for all $\phi_k\in H^{k+\frac{1}{2}}_{0,0}(\R^3_+),\ k=0,1,$ we can find 
\[u[\phi]\in L^2(\R^+_t;H^{2}_0(\R_+^3)))\cap H^2(\R^+_t;L^2(\R_+^3)))\]
such that $(u[\phi]_{\vert t=0},\partial_tu[\phi]_{\vert t=0})=\phi=(\phi_1,\phi_2)$ and 
\begin{equation}\label{l1h}\begin{array}{l}\norm{u[\phi]}_{H^{2,2}((0,+\infty)\times\R_+^3)}
\leq C\norm{\phi}_{H^{\frac{3}{2}}_{0,0}(\R^3_+)\times H^{\frac{1}{2}}_{0,0}(\R^3_+)}.\end{array}\end{equation}
Notice that, for $j+k=1$, we have
\[\hspace{-2cm}\begin{array}{l}\norm{(x_1)^{-\frac{1}{2}}\partial_{{x}_1}^j\partial_t^k v[f]_{\vert t=0}}_{L^2(\R_+^3)}\\
\leq \norm{(x_1)^{-\frac{1}{2}}(\partial_t^k{g_j}_{\vert t={x}_1})}_{L^2(\R_+^3)}+\norm{(x_1)^{-\frac{1}{2}}(\partial_{{x}_1}^j\partial_t^kv[f]_{\vert t=0}-\partial_t^k{g_j}_{\vert t={x}_1})}_{L^2(\R_+^3)}.\end{array}\]
Clearly, for $j+k=1$, we find
\[\norm{(x_1)^{-\frac{1}{2}}(\partial_t^k{g_j}_{\vert t={x}_1})}_{L^2(\R_+^3)}\leq \norm{f}_{K_2}\]
and, from Theorem 2.2 in Chapter 4 of \cite{LM2}, we get
\[\hspace{-2cm}\begin{array}{ll}\norm{(x_1)^{-\frac{1}{2}}(\partial_{{x}_1}^j\partial_t^kv[f]_{\vert t=0}-\partial_t^k{g_j}_{\vert t={x}_1})}_{L^2(\R_+^3)}&\leq C\norm{v[f]}_{H^{2,2}((0,+\infty)\times\R_+^3)}\\
\ &\leq C\norm{f}_{H^{\frac{3}{2},\frac{3}{2}}((0,+\infty)\times\R^2)}.\end{array}\]
From these two estimates we deduce that 
\[\norm{(x_1)^{-\frac{1}{2}}\partial_{{x}_1}^j\partial_t^k v[f]_{\vert t=0}}_{L^2(\R_+^3)}\leq C \norm{f}_{K_2}.\]
In the same way we show that
\[\norm{(x_1)^{-\frac{1}{2}}\partial_{({x}_2,{x}_3)}^{\beta}v[f]_{\vert t=0}}_{L^2(\R_+^3)}\leq C \norm{f}_{K_2},\ \beta\in\mathbb N^2,\ \abs{\beta}=1.\]
Combining these estimates with Theorem 2.1 in Chapter 4 of \cite{LM2} and  (\ref{l1a}), we deduce that
\[(v[f]_{\vert t=0}, \partial_t v[f]_{\vert t=0})\in H^{\frac{3}{2}}_{0,0}(\R^3_+)\times H^{\frac{1}{2}}_{0,0}(\R^3_+)\]
and we obtain
\begin{equation}\label{l1i}\norm{(v[f]_{\vert t=0}, \partial_t v[f]_{\vert t=0})}_{H^{\frac{3}{2}}_{0,0}(\R^3_+)\times H^{\frac{1}{2}}_{0,0}(\R^3_+)}\leq C \norm{f}_{K_2}.\end{equation}
In view of (\ref{l1a}), (\ref{l1h}) and (\ref{l1i}), if we set
$w[f]=v[f]+u[\phi]$ with $\phi=(-v[f]_{\vert t=0}, -\partial_t v[f]_{\vert t=0})$, conditions (\ref{l1f}) and (\ref{l1g}) will be fulfilled. \begin{flushright}
\rule{.05in}{.05in}
\end{flushright}

\textbf{Proof of Theorem \ref{t3}.} Using  local coordinate, in the same way as in the beginning of this subsection (see also Section 7.2 and the proof of Theorem 8.3 in Chapter 1 of \cite{LM1}), we can replace the space $L$ by $K_2$. Then, we  deduce Theorem \ref{t3} from the last part of Lemma \ref{l1}. \begin{flushright}
\rule{.05in}{.05in}
\end{flushright} 

\subsection{Proof of Theorem \ref{t2}}
We will now go back to Theorem \ref{t2}. First, using Theorem \ref{t3} we split  $u$ into two terms $u=v[F]+w$ with $w\in H^{2,2}(Q)$ satisfying (\ref{t3a}), (\ref{t3b}) and $v[F]$ solution of
\begin{equation}\label{eq2}\left\{\begin{array}{ll}\partial_t^2v[F]-\Delta v[F]+qv[F]=F,\quad &t\in(0,T),\ x'\in\omega,\ x_3\in\R,\\  v[F](0,\cdot)=0,\quad \partial_tv[F](0,\cdot)=0,\quad &\textrm{in}\ \Omega,\\ v[F]=0,\quad &\textrm{on}\ \Sigma\end{array}\right.\end{equation}
with $F=-(\partial_t^2w-\Delta w+qw)\in L^2(Q)$. In view of Theorem 8.1 and 8.3 in Chapter 3 of [LM1] problem (\ref{eq2}) admits a unique solution $v\in\mathcal C([0,T];H^1_0(\Omega))\cap \mathcal C^1([0,T];L^2(\Omega))$ satisfying
\begin{equation}\label{t2c}\norm{v}_{\mathcal C([0,T];H^1_0(\Omega))}+\norm{v}_{\mathcal C^1([0,T];L^2(\Omega))}\leq C\norm{F}_{L^2(Q)}\leq C\norm{w}_{H^{2,2}(Q)}.\end{equation}
\begin{remark}\label{r1}
It is well known  that  the Laplacian in $\Omega$ with Dirichlet boundary condition is a self adjoint operator associated to the sesquilinear closed coercive  form
$b$ with domain $D(b)=H^1_0(\Omega)$ and defined by 
\[b(f,g)=\int_\Omega \nabla f\cdot\overline{\nabla g}\d x.\]
Therefore, we can apply the theory introduced in Section 8 of Chapter 3 of [LM1], by considering the sesquilinear form
\[a(u,v)=\int_\Omega(\nabla_xu\cdot\overline{\nabla_xv}+qu\overline{v})\d x\]
and the space $V=H^1_0(\Omega)$, $H=L^2(\Omega)$.
\end{remark}

In view of Theorem 3.1 in Chapter 1 of [LM1] , we have\[\hspace{-2,5cm}\begin{array}{ll}w\in H^{2,2}(Q)&\subset \mathcal C([0,T];H^{\frac{3}{2}}(\Omega))\cap \mathcal C^1([0,T];H^{\frac{1}{2}}(\Omega))\\
\ &\subset \mathcal C([0,T];H^1(\Omega))\cap \mathcal C^1([0,T];L^2(\Omega))\end{array}\] and
\[\hspace{-2,5cm}\begin{array}{ll}\norm{w}_{\mathcal C([0,T];H^1(\Omega))}+\norm{w}_{\mathcal C^1([0,T];L^2(\Omega))}&\leq C(\norm{w}_{\mathcal C([0,T];H^{\frac{3}{2}}(\Omega))}+\norm{w}_{\mathcal C^1([0,T];H^{\frac{1}{2}}(\Omega))})\\
\ &\leq C\norm{w}_{H^{2,2}(Q)}.\end{array}\]
Therefore, $u\in \mathcal C([0,T];H^1(\Omega))\cap \mathcal C^1([0,T];L^2(\Omega))$ and (\ref{t3b}) implies
\[\norm{u}_{\mathcal C([0,T];H^1(\Omega))}+\norm{u}_{\mathcal C^1([0,T];L^2(\Omega))}\leq C\norm{f}_{L}.\]

It remains to show that $\partial_\nu u\in L^2(\Sigma)$ and $\norm{\partial_\nu u}_{L^2(\Sigma)}\leq C\norm{f}_L$. For this purpose, notice that $\partial_\nu u =\partial_\nu v[F]$. Thus, in view of (\ref{t3b}), the proof will be complete if we show that for $v[F]$ solution of (\ref{eq2}) we have $\partial_\nu v[F]\in L^2(\Sigma)$ and
\begin{equation}\label{t2d}\norm{\partial_\nu v[F]}_{L^2(\Sigma)}\leq C\norm{F}_{L^2(Q)}. \end{equation}
In order to prove this result, we will first assume that $q=0$ and $F$ is smooth so that $v[F]$ is sufficiently smooth. Then we will conclude by density and by repeating the arguments of Theorem A.2 in \cite{BCY}. Without lost of generality, we can also assume that $v=v[F]$ is real valued.

Let $\gamma_1\in \mathcal C^\infty(\omega,\R^2)$ be such that $\gamma_1=\nu_1$ on $\partial\omega$ with $\nu_1$ the unit outward
normal vector to  $\partial \omega$. We set $\gamma\in \mathcal C^\infty(\Omega,\R^3)\cap W^{\infty,\infty}(\Omega,\R^3)$ defined by $\gamma(x',x_3)=(\gamma_1(x'),0)$, $x'\in\omega$, $x_3\in\R$ and we obtain $\gamma=\nu$ on $\partial\Omega$. We have

\begin{equation}\label{t2e}\int_Q (\partial_t^2v-\Delta v)\gamma\cdot\nabla v\d x\d t=0.\end{equation}
Integrating by parts in $t$ we get
\[\hspace{-2,5cm}\begin{array}{ll}\int_0^T\int_\Omega \partial_t^2v\gamma\cdot\nabla v\d x\d t&=\int_\Omega \partial_tv(T,x)\gamma\cdot\nabla v(T,x)\d x-\int_0^T\int_\Omega \partial_tv\gamma\cdot\nabla \partial_tv\d x\d t\\
\ \\
\ &=\int_\Omega \partial_tv(T,x)\gamma\cdot\nabla v(T,x)\d x-\frac{1}{2}\int_0^T\int_\Omega \gamma\cdot\nabla (\partial_tv)^2\d x\d t.\end{array}\]
Now notice that $\gamma\cdot\nabla (\partial_tv)^2=\gamma_1\cdot\nabla_{x'} (\partial_tv)^2$. Therefore, applying the Green formula in $x'\in\omega$ we get
\[\hspace{-2,5cm}\begin{array}{ll}\int_0^T\int_\Omega \gamma\cdot\nabla v_t^2\d x\d t&=\int_0^T\int_\R\int_\omega \gamma_1\cdot\nabla_{x'} (\partial_tv)^2\d x'\d x_3\d t\\
\ \\
\ &=\int_0^T\int_\R\int_{\partial\omega}  v_t^2\d x'\d x_3\d t-\int_0^T\int_\R\int_\omega \textrm{div}(\gamma) v_t^2\d x'\d x_3\d t\\
\ \\
\ &=\int_{\Sigma}v_t^2\d\sigma(x)\d t-\int_Q\textrm{div}(\gamma) v_t^2\d x\d t.\end{array}\]
Since $v_{\vert\partial\Omega}=0$ we deduce that $\partial_tv_{\vert\partial\Omega}=0$ and it follows
\[\int_{\Sigma}v_t^2\d\sigma(x)\d t=0.\]
Thus we have
\begin{equation}\label{t2f}\hspace{-2,5cm}\int_0^T\int_\Omega \partial_t^2v\gamma\cdot\nabla v\d x\d t=\int_\Omega \partial_tv(T,x)\gamma\cdot\nabla v(T,x)\d x+\frac{1}{2}\int_Q\textrm{div}(\gamma) v_t^2\d x\d t.\end{equation}
On the other hand, we get
\[\hspace{-2,5cm}\int_Q -\Delta v\gamma\cdot\nabla v\d x\d t=\int_Q -\Delta_{x'}v\gamma\cdot\nabla v\d x\d t+\int_Q-\partial_{{x}_3}^2v\gamma\cdot\nabla v\d x\d t.\]
Applying the Green formula in $x'\in\omega$, we find
\[\hspace{-2,5cm}\begin{array}{l}\int_Q -\Delta_{x'}v\gamma\cdot\nabla v\d x\d t\\
\ \\
=\int_0^T\int_\R\int_\omega -\Delta_{x'} v\gamma\cdot\nabla v\d x'\d x_3\d t\\
\ \\
=-\int_0^T\int_\R\int_{\partial\omega}\abs{\partial_\nu v}^2\d \sigma(x')\d x_3\d t +\int_0^T\int_\R\int_\omega \nabla_{x'} v\cdot\nabla_{x'}(\gamma\cdot\nabla v)\d x'\d x_3\d t\\
\ \\
=-\int_{\Sigma}\abs{\partial_\nu v}^2\d\sigma(x)\d t+\int_Q\nabla_{x'} v\cdot\nabla_{x'}(\gamma\cdot\nabla v)\d x\d t.\end{array}\]
and integrating by parts in $x_3\in\R$ we obtain
\[\hspace{-2,5cm}\int_Q-\partial_{{x}_3}^2v\gamma\cdot\nabla v\d x\d t=\int_Q\partial_{{x}_3}v\partial_{{x}_3}(\gamma\cdot\nabla v)\d x\d t.\]
It follows
\[\hspace{-2,5cm}\int_Q -\Delta v\gamma\cdot\nabla v\d x\d t=-\int_{\Sigma}\abs{\partial_\nu v}^2\d\sigma(x)\d t+\int_Q\nabla_x v\cdot\nabla_x(\gamma\cdot\nabla_{x} v)\d x\d t.\]
Recall that
\[\hspace{-2,5cm}\nabla v\cdot\nabla (\gamma\cdot\nabla v)=(H\nabla v)\cdot\nabla v+\frac{1}{2}\gamma\cdot\nabla(\abs{\nabla v}^2)\]
 where $H=(\partial_{{x}_j}\gamma_i)_{1\leq i,j\leq 3}$ and $\gamma=(\gamma_1,\gamma_2,\gamma_3)$. Therefore we have
\[\hspace{-2,5cm}\begin{array}{ll}\int_Q -\Delta v\gamma\cdot\nabla v\d x\d t=&-\int_{\Sigma}\abs{\partial_\nu v}^2\d\sigma(x)\d t+\int_Q(H\nabla v)\cdot\nabla v\d x\d t\\
\ \\
\ &+\frac{1}{2}\int_Q\gamma\cdot\nabla(\abs{\nabla v}^2)\d x\d t.\end{array}\]
The Green formula in $x'\in\omega$ implies
\[\hspace{-2,5cm}\begin{array}{ll}\int_\omega\gamma\cdot\nabla(\abs{\nabla v}^2)\d x'=\int_\omega\gamma_1\cdot\nabla_{x'}(\abs{\nabla v}^2)\d x'
=\int_{\partial\omega}\abs{\nabla v}^2\d \sigma(x')-\int_\omega\textrm{div}(\gamma)\abs{\nabla v}^2\d x'.\end{array}\]
Since $v_{\vert \partial\Omega}=0$, we obtain $\abs{\nabla v}^2=\abs{\partial_\nu v}^2$ on $\Sigma$. It follows
\[\hspace{-2,5cm}\int_Q\gamma\cdot\nabla(\abs{\nabla v}^2)\d x\d t=\int_{\Sigma}\abs{\partial_\nu v}^2\d\sigma(x)\d t-\int_Q\textrm{div}(\gamma)\abs{\nabla v}^2\d x\d t\]
and we get
\begin{equation}\label{t2g}\hspace{-2,5cm}\begin{array}{ll}\int_Q -\Delta v\gamma\cdot\nabla v\d x\d t=&-\frac{1}{2}\int_{\Sigma}\abs{\partial_\nu v}^2\d\sigma(x)\d t+\int_Q(H\nabla v)\cdot\nabla v\d x\d t\\
\ &-\frac{1}{2}\int_Q\textrm{div}(\gamma)\abs{\nabla v}^2\d x\d t.\end{array}\end{equation}
Combining (\ref{t2e}), (\ref{t2f}) and (\ref{t2g}) we deduce that
\[ \hspace{-2,5cm}\begin{array}{ll}\int_{\Sigma}\abs{\partial_\nu v}^2\d\sigma(x)\d t=&2\int_Q(H\nabla v)\cdot\nabla v\d x\d t-\int_Q\textrm{div}(\gamma)\abs{\nabla v}^2\d x\d t\\
\ &+2\int_\Omega \partial_tv(T,x)\gamma\cdot\nabla v(T,x)\d x+\int_Q\textrm{div}(\gamma) v_t^2\d x\d t\end{array}\]
 and it follows
 \[\norm{\partial_\nu v}_{L^2(\Sigma))}\leq C(\norm{v}_{\mathcal C([0,T];H^1_0(\Omega))}+\norm{v}_{\mathcal C^1([0,T];L^2(\Omega))})\leq C\norm{F}_{L^2(Q)}.\]
By density we can extend this result to $F\in L^2(Q)$ and, in view of Theorem A.2 in \cite{BCY}, we deduce that this result holds for $q\neq0$. 
 \begin{flushright}
\rule{.05in}{.05in}
\end{flushright}

%%%%%%%%%%%%%%%%%%%%%%%%%%%%%%%%
\section*{References}

%%%%%%%%%%%%%%%%%%%%%%%%%%%%%%%%%%
   
%

\end{document}